\documentclass[preprint,12pt]{elsarticle}

\usepackage{amssymb}
\usepackage{amsmath}

\journal{Nuclear Physics B}
\begin{document}

\begin{frontmatter}



\title{Heged\"{u}s Szil\`agyi Fixed Point Theorem by  Hardy Rogers  and  Interpolative approach}


\cortext[cor1]{Corresponding author}
\author[label1]{Anuradha Gupta}
\affiliation[label1]{organization={Department of mathematics, Delhi College of Arts and Commerce},
	addressline={University of Delhi},
	city={New Delhi},
	postcode={110023},
	state={Delhi},
	country={India}}
\ead{dishna2@yahoo.in}
\author[label2]{Rahul Mansotra\corref{cor1}}
\affiliation[label2]{organization={ Department of mathematics, Faculty of Mathematical sciences},
	addressline={University of Delhi},
	city={New Delhi},
	postcode={110007},
	state={Delhi},
	country={India}}
\ead{mansotrarahul2@gmail.com}

\begin{abstract}
In this article, three new types (I), (II), (III) of Heged\"{u}s Szil\`agyi contraction on metric spaces by  interpolative and Hardy-Rogers methods have been introduced. Moreover, we give the fixed point theorems in this setting,  along with examples to justify the novelty of these contractions. 
\end{abstract}




\end{frontmatter}



\section{Introduction and Preliminaries}
 Thoroughout the paper, $(\mathcal{X},\mathcal{D})$ is considered   a metric space and  $f$ be a self mapping on  $\mathcal{X}$.\\
\textbf{Notation}:\begin{itemize}
	\item$\mathbb{N}$ be the set of positive integers, \item$O_{f}(x)$ denotes the iterates of $f$  at $x$,
	\item $Fix(f)$ be the collection of fixed points of $f$,  
	\item 
	$\mathcal{D}_{f}(x):= \sup \{\mathcal{D}(u,v): u,v \in O_{f}(x)\}$,\item
	$\mathcal{D}_{f}(x,y):= \sup \{\mathcal{D}(u,v): u,v \in O_{f}(x) \cup O_{f}(y) \}$,
	\item $M_{f}(x,y)=\dfrac{\mathcal{D}_{f}(x)+\mathcal{D}_{f}(y)}{2},$
	\item $\Phi$ is the collection of all self mappings $\phi$ on $[0.\infty)$ having the following properties: \begin{enumerate} \item $\phi(t)<t,$ for each $t>0$;\item for each $\epsilon >0,$ there exists $\mathfrak{d}>0$ such that for any $t \in (\epsilon, \epsilon + \mathfrak{d})$ implies $\phi(t)\leq \epsilon$  .\end{enumerate}\end{itemize}
\noindent Here we  assume $ \mathcal{D}_{f}(x)$ is finite for each $x \in \mathcal{X}.$ Further, a metric space $(\mathcal{X},\mathcal{D})$ is said to be complete if every Cauchy sequence in $\mathcal{X}$ converges in $(\mathcal{X},\mathcal{D})$.\\
\noindent The Banach Contraction Principle \cite{bm1} is the fundamental result in the Metric Fixed point theory. As Banach contraction is a continuous map, therefore, many authors try to find contractions which are not continuous. Following that path,  Hardy and  Rodgers \cite{bm2} in 1973  introduced the  most general contraction condition combining  various  contractions which need not be continuous and also gave a fixed point theorem which is given below:\\
\textbf{Theorem 1.} \cite{bm2} Let    $f$ be a   self map on a complete metric space $(\mathcal{X},\mathcal{D})$ which  satisfies Hardy-Rodgers contraction i.e. \vspace{-0.2cm}
  $$\mathcal{D}(fx,fy)\leq \alpha \mathcal{D}(x,y)+ \beta \mathcal{D}(x,fx)+ \gamma \mathcal{D}(y,fy) + \delta \mathcal{D}(x,fy)+\mu \mathcal{D}(fx,y) \vspace{-0.2cm}$$ for all $x,y \in \mathcal{X}$, where $0 \leq \alpha,\beta,\gamma,\delta,\mu$ such that $\alpha+\beta+\gamma+\delta+\mu<1$. Then, $f$ has a unique fixed  point.\\
  In 1980, Heged\"{u}s and  Szil\`agyi \cite{bm3} introduced  a new type of contraction in which they had taken diameter of orbits of corresponding points and proved fixed point theorem which is given as:\\
  \textbf{Theorem 2.} \cite{bm3} Let $f$ be a self mapping on a complete metric space $(\mathcal{X},\mathcal{D})$. Suppose there exists a function $\phi : [0,\infty)\to0,\infty)$ such that   the mapping $f$ satisfies the following:
  \begin{itemize}
  	\item[(i)] $\phi(t)<t$ for each $t>0,$
  	 \item[(ii)] $\phi$ is upper-semicontinous from the right,
  	 \item[(iii)] $\mathcal{D}(fx,fy)\leq \phi(\mathcal{D}_{f}(x,y)),$ for all $x,y \in \mathcal{X}.$
  	\end{itemize} Then, $f$ has a unique fixed point. In fact, fixed point  is approximated by a sequence of iterates of $f$ at any point.\\
 Following this Suzuki \cite{bm6} modified the result by replacing the  function $\phi$  by another
  function in $\Phi$. Consequently, Touail and Moutawakil \cite{bm7} used  generalized contractive condition   and gave  fixed point theorem which  was given as follows: \\
  \textbf{Theorem 3. } \cite{bm7} Assume that  there exists $\phi \in \Phi$ such that the self mapping $f$ on a complete metric space $(\mathcal{X},\mathcal{D})$ satisfies the following condition: $$\mathcal{D}(fx,fy)\leq \max \{ \phi(\mathcal{D}_{f}(x)),\phi(\mathcal{D}_{f}(y)),\phi(\mathcal{D}_{f}(x,y))\}$$ for all $x,y \in \mathcal{X},\phi \in \Phi$.  Then, $f$ has a unique fixed point.
  
Recently,  interpolative contraction was introduced by Karapinar \cite{bm4} on metric spaces. Interpolative  contractions are the  most refined version  of older   contractions which need not have a unique fixed point. In 2020, Karapinar et.al \cite{bm5}  introduced Interpolative Boyd-Wong and Matkowski type contractions on metric spaces and provided fixed point theorem  which was given in the following:\\
\textbf{Theorem 4.} \cite{bm5} If a self-mapping $f$   on a complete metric space $(\mathcal{X},\mathcal{D})$ satisfies the interpolative Boyd-Wong contraction i.e. $$ \mathcal{D}(fx,fy) \leq \phi\Big(\mathcal{D}^{\alpha}(x,y)\mathcal{D}^{\beta}(x,fx)\mathcal{D}^{\gamma}(y,fy)\Big[\dfrac{\mathcal{D}(fx,y)+\mathcal{D}(x,fy)}{2}\Big]^{1-\lambda}\Big) $$ or interpolative Matkowski contraction i.e.
 $$ \mathcal{D}(fx,fy) \leq \psi\Big(\mathcal{D}^{\alpha}(x,y)\mathcal{D}^{\beta}(x,fx)\mathcal{D}^{\gamma}(y,fy)\Big[\dfrac{\mathcal{D}(fx,y)+\mathcal{D}(x,fy)}{2}\Big]^{1-\lambda}\Big),  $$ where $\alpha,\beta,\gamma >0$ and $\lambda=\alpha+\beta+\gamma<1$, $\phi $ and $\psi$ be two self maps on $[0,\infty)$ such that $\phi$ and $\psi$ are increasing, $\phi(t)<t$, and $\lim \limits_{n \to \infty}\psi^{n}(t)=0$,  for each $t>0$, $\phi$ is upper-semicontinuous. Then, $f$ has atleast one fixed point.\\
  Motivated by the works of Hardy  and Rodgers,  Karapinar et al., and Touail and Moutawakil on  contractions, we introduce three new types of contractions on  metric spaces. Fixed point theorems are also provided  with examples in support of the results obtained.
\section{Main Results}
We begin this section with the following definition:
 
\noindent\textbf{Definition 5}. A self map $f$ on $\mathcal{X}$  defines a generalized contraction of type (I) if there exists $\phi \in \Phi,$ and   $0<\alpha, \beta, \gamma $, and $ \lambda=\alpha+\beta+\gamma <1 $ such that   
 \begin{equation*}\mathcal{D}(fx,fy)\leq \phi( \mathcal{D}_{f}(x))^{\alpha} \phi (\mathcal{D}_{f}(y) )^{\beta}\phi (\mathcal{D}_{f}(x,y))^{\gamma} \phi(M_{f}(x,y))^{1-\lambda},\end{equation*} for all $x,y \in \mathcal{X}\textbackslash Fix(f)$.\\
   \textbf{Theorem 6}. Let $f$ be an  a generalized contraction of type (I) on a complete metric space $(\mathcal{X},\mathcal{D})$. Then it  has atleast one  fixed point.\\
   \textbf{Proof.} Let $x \in \mathcal{X}$  and consider the sequence $(x_{n})$ defined as $f^{n}x=x_{n}$, for each $n \in \mathbb{N}$. If $x_{n}=x_{n+1}$ for some $n \in \mathbb{N}$, then $fx_{n}=x_{n}$. Assume $x_{n}\neq x_{n+1}$ for all $n \in \mathbb{N}$. Also, $  O_{f}(x_{n+1}) \subseteq  O_{f}(x_{n})$ for all $n \in \mathbb{N}.$ Thus, $\mathcal{D}_{f}(x_{n+1}) \leq \mathcal{D}_{f}(x_{n})$, for all $n \in \mathbb{N}.$ Hence, $\lim \limits_{n\to \infty}\mathcal{D}_{f}(x_{n})=\epsilon\geq 0.$ Keeping in mind $ M_{f}(x,y)=\dfrac{\mathcal{D}_{f}(x)+\mathcal{D}_{f}(y)}{2}$, for each $x,y\in \mathcal{X}$. We claim $\epsilon=0$. On the contrary suppose $\epsilon>0.$ Then, either $\epsilon < \mathcal{D}_{f}(x_{n})$, for each $n \in \mathbb{N}$ or $\epsilon = \mathcal{D}_{f}(x_{n})$ for some $n \in \mathbb{N}.$ \\In  the first case, we can take  some $k \in \mathbb{N}$ and $\mathfrak{d} >0$ such that $\mathcal{D}_{f}(x_{k})< \epsilon + \mathfrak{d}.$ Let $m,n \geq k$, for $m,n \in \mathbb{N}$. Then, we have  $$\epsilon <\mathcal{D}_{f}(x_{\max\{m,n\}})  \leq \mathcal{D}_{f}(x_{m},x_{n}) \leq \mathcal{D}_{f}(x_{\min\{m,n\}}) \leq \mathcal{D}_{f}(x_{k})< \epsilon + \mathfrak{d},$$ hence by the  property (ii) of $\phi$, we get  \begin{align*}\mathcal{D}(fx_{m},fx_{n})&\leq \phi( \mathcal{D}_{f}(x_{m}))^{\alpha} \phi (\mathcal{D}_{f}(x_{n}) )^{\beta}\phi (\mathcal{D}_{f}(x_{m},x_{n}))^{\gamma} \phi(M_{f}(x_{m},x_{n}))^{1-\lambda} \\&  \leq \epsilon^{\alpha+\beta+\gamma+1-\lambda}=\epsilon,\end{align*}
   which is  contrary to $\epsilon <  \mathcal{D}_{f}(x_{k+1})$.\\ On the other hand, if $\epsilon = \mathcal{D}_{f}(x_{k})$  holds for some $k \in \mathbb{N}$. Let $m,n \geq k$, for $m,n \in \mathbb{N}$. Then, $$\epsilon \leq\mathcal{D}_{f}(x_{\max\{m,n\}})  \leq \mathcal{D}_{f}(x_{m},x_{n}) \leq \mathcal{D}_{f}(x_{\min\{m,n\}}) \leq \mathcal{D}_{f}(x_{k})= \epsilon,$$ implies by the property  of $\phi$, we get 
   \begin{align*}\mathcal{D}(fx_{m},fx_{n})&\leq \phi( \mathcal{D}_{f}(x_{m}))^{\alpha} \phi (\mathcal{D}_{f}(x_{n}) )^{\beta}\phi (\mathcal{D}_{f}(x_{m},x_{n}))^{\gamma} \phi(M_{f}(x_{m},x_{n}))^{1-\lambda}\\ &= \phi( \epsilon)^{\alpha+\beta+\gamma+ 1-\lambda} <\epsilon,    \end{align*} which is a contrary to $\epsilon <  \mathcal{D}_{f}(x_{k+1})$. Therefore, $\lim \limits_{n \to \infty}\mathcal{D}_{f}(x_{n})=0$.\\ Similarly, let if possible,  the sequence $( f^{n}y=y_{n})$ with $y_{n}\neq y_{n+1} $, for all $n \in \mathbb{N}$ and $y \in \mathcal{X}$ such that $\lim \limits_{n \to \infty}\mathcal{D}_{f}(y_{n})=0$. Now we will prove that $\lim \limits_{n \to \infty}\mathcal{D}_{f}(x_{n},y_{n})=0.$  Then  $O_{f}(x_{{n+1}}) \cup O_{f}(y_{{n+1}}) \subseteq  O_{f}(x_{n}) \cup O_{f}(y_{n})$, for each $n \in \mathbb{N}$. Thus, $\mathcal{D}_{f}(x_{n+1},y_{n+1})\leq\mathcal{D}_{f}(x_{n},y_{n})$, for all $n \in \mathbb{N}$. So, $\lim \limits_{n \to \infty}\mathcal{D}_{f}(x_{n},y_{n})=\epsilon$. If $\epsilon >0$, then atleast one of the following holds: $\epsilon<\mathcal{D}_{f}(x_{n},y_{n})$ for all $n \in \mathbb{N}$ or $\epsilon =\mathcal{D}_{f}(x_{n},y_{n}) $ for some $n \in \mathbb{N}.$ Suppose $\epsilon<\mathcal{D}_{f}(x_{n},y_{n})$ for all $n \in \mathbb{N}$ holds, then we can select $\mathfrak{d}>0$ and  $k \in \mathbb{N}$ such that $\mathcal{D}_{f}(x_{k},y_{k}) < \epsilon +\mathfrak{d}$ and $\mathcal{D}_{f}(x_{k})< \epsilon \text{ and }\mathcal{D}_{f}(y_{k})< \epsilon$ as  $\lim \limits_{n \to \infty}\mathcal{D}_{f}(x_{n})=0 $ and  $\lim \limits_{n \to \infty}\mathcal{D}_{f}(y_{n})=0 $. For $m,n\geq k$, we have \vspace{-0.2cm} \begin{gather*}\epsilon <  \mathcal{D}_{f}(x_{\max\{m,n\}},y_{\max\{m,n\}}) \leq  \mathcal{D}_{f}(x_{m},y_{n}) \leq \mathcal{D}_{f}(x_{\min\{m,n\}},y_{\min\{m,n\}})\\\leq \mathcal{D}_{f}(x_{k},y_{k})<\epsilon +\mathfrak{d},  \vspace{-0.5cm}\end{gather*} and by the property (ii) of $\phi$, we get \begin{align*}\mathcal{D}(fx_{{m}},fy_{n})&\leq \phi( \mathcal{D}_{f}(x_{m}))^{\alpha} \phi (\mathcal{D}_{f}(y_{n}) )^{\beta}\phi (\mathcal{D}_{f}(x_{m},y_{n}))^{\gamma} \phi(M_{f}(x_{m},y_{n}))^{1-\lambda}\\& <  \mathcal{D}_{f}(x_{m})^{\alpha}  \mathcal{D}_{f}(y_{n} )^{\beta}\epsilon^{\gamma} M_{f}(x_{m},y_{n})^{1-\lambda}\\&<\epsilon^{\alpha+\beta+\gamma+1-\lambda}=\epsilon,\vspace{-0.3cm}\end{align*}
   which is  contrary to the fact $\epsilon<\mathcal{D}(x_{{k+1}},y_{k+1}).$ On the other side if $\epsilon =\mathcal{D}(x_{{k}},y_{k})\text{ and }\mathcal{D}_{f}(y_{k})< \epsilon$, for some $k\in \mathbb{N}$ as  $\lim \limits_{n \to \infty}\mathcal{D}_{f}(x_{n})=0 $ and  $\lim \limits_{n \to \infty}\mathcal{D}_{f}(y_{n})=0 $, then for $m,n\geq k,$ we have \begin{gather*} \epsilon \leq  \mathcal{D}_{f}(x_{\max\{m,n\}},y_{\max\{m,n\}}) \leq  \mathcal{D}_{f}(x_{m},y_{n}) \leq \mathcal{D}_{f}(x_{\min\{m,n\}},y_{\min\{m,n\}})\\\leq \mathcal{D}_{f}(x_{k},y_{k})=\epsilon,  \end{gather*}
    and by the property of $\phi$, we get \begin{align*}\mathcal{D}(fx_{{m}},fy_{n})&\leq \phi( \mathcal{D}_{f}(x_{m}))^{\alpha} \phi (\mathcal{D}_{f}(y_{n}) )^{\beta}\phi (\mathcal{D}_{f}(x_{m},y_{n}))^{\gamma} \phi(M_{f}(x_{m},y_{n}))^{1-\lambda}\\& <  \mathcal{D}_{f}(x_{m})^{\alpha}  \mathcal{D}_{f}(y_{n} )^{\beta}\phi (\epsilon)^{\gamma} M_{f}(x_{m},y_{n})^{1-\lambda}\\&<\epsilon^{\alpha+\beta+\delta+1-\lambda}=\epsilon,\end{align*}
   which is  contrary to  $\epsilon=\mathcal{D}(x_{{k+1}},y_{k+1}).$ Hence, $ \lim \limits_{n \to \infty }\mathcal{D}_{f}(x_{n},y_{n})=0$.
   Next we will show the existence of fixed point.  Since $\lim \limits_{n \to \infty }\mathcal{D}_{f}(x_{n})=0$, then $(x_{n})$ is a Cauchy sequence in $\mathcal{X}$. As $(\mathcal{X},\mathcal{D})$ is complete, so there exists $w \in \mathcal{X}$ such that $\lim \limits_{n \to \infty}x_{n}=w.$ On the similar lines if possible we can consider the sequence $(w_{n}=f^{n}w)$ such that $w_{n}\neq w_{n+1}$, for all $n \in \mathbb{N}$, then we have $\lim \limits_{n \to \infty}\mathcal{D}_{f}(x_{n},w_{n})=0$, which gives $\lim \limits_{n \to \infty}w_{n}=w.$ We  claim that $\mathcal{D}_{f}(w)=0.$      Suppose $\mathcal{D}_{f}(w)=\epsilon>0. $ As $\lim \limits_{n \to \infty}\mathcal{D}_{f}(w_{n})=0, $ so there is $k\in \mathbb{N}$ such that $
   	\epsilon =\mathcal{D}_{f}(w)=...=\mathcal{D}_{f}(w_{k})> \mathcal{D}_{f}(w_{k+1}).$ Now, for $n>k$, we consider the following:
   	If  $n-1=k$, then $\phi(\mathcal{D}_{f}(w_{n-1}))=\phi(\mathcal{D}_{f}(w_{k}))=\phi(\epsilon).$ On the other hand, if $n-1>k$, then $\phi(\mathcal{D}_{f}(w_{n-1}))< \mathcal{D}_{f}(w_{n-1})\leq \mathcal{D}_{f}(w_{k+1})<\epsilon.$ Further, 
   	 \begin{gather*}\mathcal{D}(fw_{k-1},fw_{n-1})\leq \phi( \mathcal{D}_{f}(w_{k-1}))^{\alpha}  \phi (\mathcal{D}_{f}(w_{n-1}) )^{\beta} \phi (\mathcal{D}_{f}(w_{k-1},w_{n-1}))^{\gamma}\\\phi(M_{f}(w_{k-1},w_{n-1}))^{1- \lambda}< \phi( \mathcal{D}_{f}(w_{k-1}))^{\alpha}  \phi (\mathcal{D}_{f}(w_{n-1}) )^{\beta} \phi (\mathcal{D}_{f}(w_{k-1}))^{\gamma}\\M_{f}(w_{k-1},w_{n-1})^{1- \lambda}< \phi(\epsilon)^{\alpha+\gamma}\epsilon^{\beta+1-\lambda}<\epsilon. \end{gather*}
   	Thus, $\epsilon = \sup \{ \mathcal{D}(w_{k},w_{n}): n>k \}< \epsilon$ which is a contradiction.	Hence $\mathcal{D}_{f}(w)=0.$ Moreover, if the  generalized contraction of type (I) holds on the fixed point set of $f$, then fixed point is unique. \\
 \textbf{Example 7}. Let $\mathcal{X}=[0,\infty)$, and $a,b>0$, and $\mathcal{D}$ denotes the usual metric on $\mathcal{X}$. Now we define 
   $f$ as a self map on $\mathcal{X}$  as:
 \begin{equation*}
 fx=\begin{cases}b, \text{ if } x \neq a, \\
 a, \text{ else. } \end{cases}
 \end{equation*}
 and the self map $\phi $ on $[0,\infty)$  as $\phi(t)= \dfrac{t}{2}$, for all $t \in [0,\infty)$. Clearly, $\phi(t)<t$,   for each $t>0.$ Also, for any $\epsilon>0$, choose $\mathfrak{d}=\epsilon$, then $$\epsilon < t< \epsilon +\mathfrak{d}=2  \epsilon \text{ implies } \phi(t)=\dfrac{t}{2} < \epsilon. $$ So, $\phi \in \Phi$. 
 It can be easily verified that $f$ is a generalized contraction of type (I)  on $(\mathcal{X},\mathcal{D})$ for any $0<\alpha, \beta, \gamma$  as  $\mathcal{D}(fx,fy)=0,$ for each $x,y \in \mathcal{X}\textbackslash Fix(f).$  Clearly $b$ and $a$ are the fixed points of $f$. \\
  \textbf{Definition 8}. A self map $f$ on $\mathcal{X}$ defines a generalized contraction of type (II) if there exists $\phi \in \Phi,$ and   $0  \leq  \alpha, \beta, \gamma,\delta  $, and $ \alpha+\beta+\gamma+\delta \leq 1 $ such that
 \begin{equation*}\mathcal{D}(fx,fy)\leq \alpha \phi( \mathcal{D}_{f}(x))+ \beta \phi (\mathcal{D}_{f}(y) )+\gamma \phi (\mathcal{D}_{f}(x,y))+ \delta\phi(M_{f}(x,y)),\end{equation*} for all $x,y \in \mathcal{X}$. \\
 \textbf{Theorem 9}. Let $f$ be a generalized contraction of type (II) on a complete metric space $(\mathcal{X},\mathcal{D})$. Then, it  has a unique fixed point. Moreover, for each $x \in \mathcal{X},$ the sequence $(f^{n}x)$ converges to the fixed point of $f$.  \\
 \textbf{Proof.} First we prove the uniqueness of fixed point of $f$. For  $x,y \in Fix(f)$,  we have  
  $\mathcal{D}(x,y)=\mathcal{D}(fx,fy)\leq \alpha \phi( \mathcal{D}_{f}(x))+ \beta \phi (\mathcal{D}_{f}(y) )+\gamma \phi (\mathcal{D}_{f}(x,y))+ \delta\phi(M_{f}(x,y))= \gamma\phi(\mathcal{D}(x,y)). $
  Hence, $x=y$. The rest of the proof follows easily from the proof of  Theorem 1.\\
  \textbf{Example 10.} 
Consider $\mathcal{X}=\mathbb{N}\backslash \{3\}$ with the usual metric $\mathcal{D}$ on it. Then $\mathcal{X}=A\cup B \cup C$, where $A=\{7,11,15,...\}, B=\{ 1,5,9,...\}, C=\{2,4,6,...\} $.   Now we define 
$f$ as a self map on $\mathcal{X}$  as:
\begin{equation*}
	fx=\begin{cases}2, \text{ if } x \in A, \\
		1, \text{ else. } \end{cases}
\end{equation*}
  and the self map $\phi $ on $[0,\infty)$   as: \begin{equation*}
  	\phi(t)=\begin{cases}0, \text{ if } t  \text{ is even, } \\
  		\dfrac{5t}{6}, \text{ else. } \end{cases}
  \end{equation*}
  Clearly, $\phi(t)<t$,   for each $t>0.$ Also, for any $\epsilon>0$, choose $\mathfrak{d}=\dfrac{\epsilon}{5}$, then $$\epsilon < t< \epsilon +\mathfrak{d}=\dfrac{6\epsilon}{5} \text{ implies } \phi(t)=0 \text{ or }\dfrac{5t}{6} < \epsilon. $$  So, $\phi \in \Phi.$ We will show that $f$ is a
 generalized contraction of type (II) for any  $0\leq \alpha,\beta,\gamma,$ and $\delta=\dfrac{5}{6}$ with $\alpha+\beta+\gamma+\delta\leq 1$ . This can be easily verified from the following table: \begin{center}
 \begin{tabular}{|c|c|}
 	\hline
  $\mathcal{D}(fx,fy)$& Lower bound of  \\&$\alpha\phi( \mathcal{D}_{f}(x))+ \beta \phi (\mathcal{D}_{f}(y) )+\gamma \phi (\mathcal{D}_{f}(x,y))+ \delta\phi(M_{f}(x,y))$ \\ \hline
 	\multicolumn{2}{|c|}{Case 1. If both  $x,y \in A $ or $\mathcal{X}\backslash A$}  \\ \hline
 	0 & 0 \\ \hline
 	\multicolumn{2}{|c|}{\hspace{-1.6cm}Case 2. If $x \in A,y =1$  } \\ \hline
 1 & $2.066$ \\ \hline
 	\multicolumn{2}{|c|}{\hspace{-1.6cm}Case 3. If $x\in A, y \in C$} \\ \hline
 $1$ & $2.41$ \\ \hline
 	\multicolumn{2}{|c|}{\hspace{-1.6cm}Case 4. If $x\in A,y\in B$} \\ \hline
 	$1$ & $2.755$ \\ \hline
 \end{tabular}\end{center}
 
 Thus, by Theorem 9, $1$ is the unique fixed point of $f$. \\
 \textbf{Remark.} Note that $\alpha+\beta+\gamma+\delta$ can be taken to be  1 as compared to Hardy-Rodgers contraction. \\
 \textbf{Definition 11.} A self mapping $f$ on a metric space $(\mathcal{X},\mathcal{D})$ defines a generalized contraction of type (III) if there exists $\phi \in \Phi$ such that $f$ satisfy the following:\vspace{-0.5cm}$$ \mathcal{D}(fx,fy)\leq \max \Big\{\phi(\mathcal{D}_{f}(x)),\phi(\mathcal{D}_{f}(y)),\phi(\mathcal{D}_{f}(x,y),\phi\Big(\dfrac{\mathcal{D}_{f}(x)+\mathcal{D}_{f}(y)}{2}\Big)\Big\}\vspace{-0.4cm}$$ for all $x,y \in \mathcal{X}.$\\
 \textbf{Theorem 12. }Let $f$ be a generalized contraction of type (III) on a complete metric space $(\mathcal{X},\mathcal{D})$. Then $f$ has a unique fixed point. Moreover, $f^{n}x$ converges to a fixed point for each $x\in \mathcal{X}.$\\
 \textbf{Proof. } For uniqueness of fixed point, let $x,y$ be fixed points of $f$. Then $\mathcal{D}(x,y)=\mathcal{D}(fx,fy)\leq \phi (\mathcal{D}_{f}(x,y))<\mathcal{D}_{f}(x,y)$. Thus, $x=y$. Rest of the proof follows from the proof of Theorem $1$.\\
 \textbf{Example 13.} The map $f$ in Example 10 is a generalized   contraction of type (II). \\
\textbf{Remark.} By Example 10 and 13, it can be seen that a generalized contractions of types (II) and  (III) are more general versions  than the  contraction in Theorem $3$ as for $x\in A,y\in B$, $\mathcal{D}(fx,fy)=1,  \phi( \mathcal{D}_{f}(x))= \phi (\mathcal{D}_{f}(y) )=\phi (\mathcal{D}_{f}(x,y))=0$.

\section{Conclusion} 
In the above  results, we ensure the existence of fixed point for generalized contractions of types (I),  (II) and  (III) on metric spaces. This methodology may be applied to other contractions as well. In addition, these concepts can be extended to more general spaces available in the literature.

\end{document}